\documentclass[12pt]{elsarticle}
\usepackage{amssymb, amsthm, enumerate, graphicx, amscd, amsmath, amssymb,
	amsfonts}

\numberwithin{equation}{section}
\newtheorem{theorem}[equation]{Theorem}

\usepackage{xcolor}

\begin{document}

\title{The critical point and the $p$-norm of $A_s$ and $C$-matrices}
	
\author{Ludovick Bouthat\corref{cor1}}
\ead{ludovick.bouthat.1@ulaval.ca}
\author{Javad Mashreghi}
\ead{javad.mashreghi@mat.ulaval.ca}

\address{D\'epartement de math\'ematiques et de statistique, Universit\'e Laval,
	Qu\'ebec, QC,
	Canada, G1K 0A6}

\cortext[cor1]{Corresponding author}


\begin{abstract}
The $L$-matrix $A_s=[1/(n+s)]$ was introduced in \cite{MRtmp}. As a surprising property, we showed that its 2-norm is constant for $s \geq s_0$, where the critical point $s_0$ is unknown but relies in the interval $(1/4,1/2)$. In this note, using some delicate calculations we sharpen this result by improving the upper and lower bounds of the interval surrounding $s_0$. Moreover, we show that the same property persists for the $p$-norm of $A_s$ matrices. We also obtain the 2-norm of a family of $C$-matrices with lacunary sequences.
\end{abstract}
	
\maketitle

\section{Introduction}
We encountered the $L$-matrices in studying the Hadamard multipliers in function spaces \cite{MR4032210}. Characterizing $\mathcal{M}(X)$,
the multiplier algebra of a Banach space $X$ of analytic functions on the open unit disc $\mathbb{D}$, is a very important subject in various studies of function spaces, e.g., zero sets, invariant subspaces, and cyclic elements \cite{MR2500010}. In \cite{MR4032210}, it is shown that $h(z) = \sum_{n=0}^{\infty} c_n z^n$ is a Hadamard multiplier for every superharmonically weighted Dirichlet Space $\mathcal{D}_\omega$ if and only if the infinite matrix
\[
T_h =
\begin{pmatrix}
c_1-c_0 & c_2-c_1 & c_3-c_2 & c_4-c_3 & \cdots \\
0 & c_2-c_1 & c_3-c_2 & c_4-c_3 & \cdots \\
0 & 0 & c_3-c_2 & c_4-c_3 & \cdots \\
0 & 0 & 0 & c_4-c_3 & \cdots \\
\vdots & \vdots & \vdots & \vdots & \ddots \\
\end{pmatrix}.
\]
acts as a bounded operator on $\ell^2$. This observation led to deeper study of $C$-matrices and $L$-matrices in \cite{MRtmp} and \cite{MRtmp2}, which are interesting subjects by themselves. More generally, it is often important to determine if an infinite matrix act as a bounded operator on some sequence space. Famous examples include the infinite Hilbert matrix \cite{MR1554854} or the Cesàro matrix \cite{MR0030620} which are important tools in approximation theory and in the study of divergent sequence, respectively.

Let $(a_n)_{n \geq 0}$ be a sequence of complex numbers. Then the infinite matrix
\[
A =
\begin{pmatrix}
a_0 & a_1 & a_2 & a_3 & \cdots \\
a_1 & a_1 & a_2 & a_3 & \cdots \\
a_2 & a_2 & a_2 & a_3 & \cdots \\
a_3 & a_3 & a_3 & a_3 & \cdots \\
\vdots & \vdots & \vdots & \vdots & \ddots \\
\end{pmatrix},
\]
is called an {\em $L$-matrix}. Abusing the notation, we will write $A=[a_n]$ and, despite being slightly confusing, a general element of $A$ will be denoted by $a_{ij}$, where $i$ and $j$ run through $\{0,1,2,\dots\}$. Note that due to connections to function space theory and Taylor series of analytic functions on $\mathbb{D}$, the indices starts from zero. See also \cite[Page 42]{MR0305568} for another class of $L$-matrices, used in the theory of large linear systems. The infinite matrix
\[
C =
\begin{pmatrix}
a_0 & 0 & 0 & 0 & \cdots \\
a_1 & a_1 & 0 & 0 & \cdots \\
a_2 & a_2 & a_2 & 0 & \cdots \\
a_3 & a_3 & a_3 & a_3 & \cdots \\
\vdots & \vdots & \vdots & \vdots & \ddots \\
\end{pmatrix},
\]
is also called a {\em $C$-matrix}, which as in the previous case will be denoted by $C[a_n]$.

To the best of our knowledge, the first instance of an $L$-matrix being used in the literature is in 1983 in an article from Choi \cite{MR701570}. In this note, the author used the special $L$-matrix $A_1$ to show that the infinite Hilbert matrix act as a bounded operator on $\ell^2$. He even suggest that this is perhaps the quickest way to show the boundedness of the infinite Hilbert matrix as an operator from $\ell^2$ to $\ell^2$.

The $C$-matrices have received considerable attention over the past 30 years. Because of the close relation with the $L$-matrices, we will keep the notation \emph{$C$-matrices} in this paper; however in most of the literature they are referred to by the name \emph{terraced matrices}. This name seems to have been introduced by Rhaly in his note from 1989 \cite{MR998639}. The same author then published six more articles on this subject, with his latest being in 2013 \cite{MR3112183}. Other examples of studies of the $C$-matrices include Roades \cite{MR1034416}, who provide lower bounds of $p$-norms of these matrices under certain restrictions on $p$, Almasri \cite{MR3458973} show that the $C$-matrix defined by the sequence $1/n^\alpha$ is $p$-summing if and only if $\alpha>1$, and Durna \& Yildirim \cite{MR3527878} introduced and studied the generalized terraced matrices.

Recently, the study of the $p$-norm of different infinite matrices has been an active research area in computational mathematics. For example, \.{I}lkhan \cite{MR4033299} study the $p$-norms of some matrix operators on Fibonacci weighted difference sequence space. Jevti\'{c} and Karapetrovi\'{c} \cite{MR3641772} obtained some result on the infinite Hilbert matrix on spaces of {B}ergman-type. A related matrix is the so-called multiplicative Hilbert matrix $M$, the infinite matrix with entries $(\sqrt{mn}\log(mn))^{-1}$, where $m,n\geq2$. In their note, Brevig and al. \cite{MR3545935} obtained some results on the $p$-norm of this matrix and even more. As a final example, we mention Chalendar and Partington \cite{MR3208804} who showed in their article that if $T$ is a bounded operator on $H^2$, then under certain natural conditions it will act as a bounded operator on $H^2(\beta)$ and it will satisfy the inequality $\|T\|_{H^2} \leq \|T\|_{H^2(\beta)}$.

\section{Main results}
In \cite{MRtmp}, we studied the $L$-matrix
\[
A_s=\left[\frac{1}{n+s}\right] =
\begin{pmatrix}
\frac{1}{s}\,   & \frac{1}{1+s}\, & \frac{1}{2+s}\, & \frac{1}{3+s}\, & \cdots\\[3pt]
\frac{1}{1+s}\, & \frac{1}{1+s}\, & \frac{1}{2+s}\, & \frac{1}{3+s}\, & \cdots\\[3pt]
\frac{1}{2+s}\, & \frac{1}{2+s}\, & \frac{1}{2+s}\, & \frac{1}{3+s}\, & \cdots\\[3pt]
\frac{1}{3+s}\, & \frac{1}{3+s}\, & \frac{1}{3+s}\, & \frac{1}{3+s}\, & \cdots\\
\vdots\, & \vdots\, & \vdots\, & \vdots\, & \ddots \\
\end{pmatrix}
\]
and, surprisingly enough, we could precisely determine the norm for some values of $s$. As a matter of fact, we showed that
\begin{equation}\label{thm: norm A_s}
\|A_s\|_{\ell^2 \to \ell^2} = 4, \qquad  \big(s \geq \frac{1}{2} \big).
\end{equation}
Since $a_{0} = 1/s$, we certainly have
\[
\|A_s\|_{\ell^2 \to \ell^2} \geq \frac{1}{s} >4, \qquad \big(0<s<\frac{1}{4}\big).
\]
Therefore, in the light of \eqref{thm: norm A_s}, an interesting question is to determine the {\em critical point} $s_0$, where
\[
s_0 := \inf \{s: \|A_s\|_{\ell^2 \to \ell^2} = 4\}.
\]
From the previous observations, we know that
\[
\frac{1}{4} \leq s_0 \leq \frac{1}{2}.
\]
We sharpen these estimations as follows.
\begin{theorem} \label{T:so-2}
We have
\[
\frac{ \sqrt{6(8+3\sqrt{3})}-\sqrt{3}-3 }{12} \leq s_0 \leq \frac{1}{2\sqrt{2}}.
\]
\end{theorem}
The above upper and lower estimations are highly non-trivial and involve delicate calculations.	
	
We also study the $p$-norm of $A_s$ defined as
\[
\|A_s\|_{\ell^p \to \ell^p} := \sup_{x \ne 0} \frac{\|A_s x\|_p}{\|x\|_p}.
\]
In this case, the same interesting phenomenon persist, albeit for $s \geq 1$.
	
\begin{theorem} \label{T:As-p}
Let $s \geq 1$. We have
\begin{align}\label{eq: p-norm A_s}
\|A_s\|_{\ell^p \to \ell^p} = \frac{p^2}{p-1}, \qquad (1<p<\infty).
\end{align}
\end{theorem}

Lastly, we study the 2-norm of a special family of $C$-matrices with {\em lacunary} coefficient.
To show that the necessary condition $a_n=O(1/\sqrt{n})$ for the $L$-matrix $A=[a_n]$ to be a bounded operator on $\ell^2$ is sharp, we introduced in \cite{MRtmp2} the $L$-matrices with lacunary coefficients. We say that the sequence $(a_n)$ is lacunary if there is a constant $\rho>1$ and a subsequence $(n_j)_{j \geq 1}$ of positive integers such that
\[
\frac{n_{j+1}}{n_{j}} \geq \rho
\]
and $a_n=0$ except possibly for indices $n \in \{n_j: j \geq 1\}$. In particular, we introduced the special Ces\`aro type matrix $C=C[a_n]$ defined by
\begin{equation}\label{E:matrix-B-lacun}
a_{N^j} = \frac{1}{N^{j/2}}, \qquad (j \geq 1),
\end{equation}
and $a_n=0$ for other values of $n$, where $N\geq2$ in a fixed integer. Then we showed that
\[
\frac{\sqrt{N}}{\sqrt{N-1} \,\,} \leq \|C\|_{\ell^2 \to \ell^2} \leq
\frac{\sqrt{N}}{\sqrt{N}-1}.
\]
This estimation is not far from being optimal. However, with more elaborate calculation, it is possible to precisely determine $\|C\|_{\ell^2 \to \ell^2}$.
	
\begin{theorem}\label{C:lacunary-Lmatrix-2}
Let $N\geq2$ be a fixed positive integer, and let $C$ be defined by \eqref{E:matrix-B-lacun}. Then
\[
\|C\|_{\ell^2 \to \ell^2} =  \frac{\sqrt{N-1}}{\sqrt{N}-1}.
\]
\end{theorem}

\section{Approximating $s_0$}
In this section, we present the proof of Theorem \ref{T:so-2}. Theorem 2 from \cite{MRtmp} tells us that if $A=[a_n]$ is an $L$-matrix and that
there is a sequence of strictly decreasing positive numbers $\delta_n$, $n \geq
0$, such that
\[
\Delta:=\sup_{n \geq 1} \frac{(|a_n|+\delta_{n-1})
	(|a_n|+\delta_{n})}{\delta_{n-1}-\delta_{n}} < \infty.
\]
Then $A \in \mathcal{L}(\ell^2)$ and, moreover,
\[
\|A\|_{\ell^2 \to \ell^2} \leq \max\{ \delta_0+|a_0|, \, \Delta \}.
\]
Consider the sequence $(\delta_n)$ defined by $\delta_n =
(n+s+\frac{1}{2})^{-1}$.

Then we know that
\begin{align*}
\|A_s\|_{\ell^2 \to \ell^2}
&\leq \max \left\{ \delta_0 + |a_0| ,\, \frac{\left( |a_n| + \delta_{n-1}
	\right) \left( |a_n| + \delta_{n} \right) }{ \delta_{n-1}-\delta_n } (n \geq 1)
\right\}\\
&=\max \left\{ \frac{1}{ s+\frac{1}{2} }+\frac{1}{s} ,\, 4- \frac{1}{ 4(n+s)^2
}(n \geq 1) \right\}\\
&=\max \left\{ \frac{1}{ s+\frac{1}{2} }+\frac{1}{s} ,\, 4 \right\}.
\end{align*}
Observe that $f(s) := \frac{1}{ s+\frac{1}{2} } + \frac{1}{s}$ is a strictly
decreasing function. Thus, if $s \geq \frac{1}{ 2\sqrt{2} }$, we have
\[
f(s) \leq f \Big( \frac{1}{ 2\sqrt{2} } \Big) = 4.
\]
Hence, we conclude that $\|A_s\|_{\ell^2 \to \ell^2} \leq 4$ whenever $s \geq
\frac{1}{ 2\sqrt{2} }$. Moreover, from the previous observations, we have
$\|A_s\|_{\ell^2 \to \ell^2} \geq 4$ whenever ${s \leq 1/2}$ so we conclude that
$s_0 \leq \frac{1}{ 2\sqrt{2} }$. It is easy to show that the method outlined
above is optimal for sequences $(\delta_n)$ of the form $\delta_n =
\alpha/(n+\beta)$.

Finding a good lower bound for the number $s_0$ turned out to be harder than the
upper bound, mainly due to the fact that we do not have a similar result to
Theorem 2 from \cite{MRtmp} for lower bounds. However, we can show that
\[
s_0 \geq \frac{\sqrt{6(8+3\sqrt{3})}-\sqrt{3}-3}{12} =: s^* \approx 0.347
\]
by using the fact that $\|A_s\|_{\ell^2 \to \ell^2} \geq \|A_s x\|_2/\|x\|_2$ for
every $x\in\ell^2$ and by carefully choosing the entries of the sequence $x$.
Note that the upper bound is $s_0 \leq 1/(2\sqrt{2}) \approx 0.354$.

Since for $s<\frac{1}{4}$, we have $\|A_s\|_{\ell^2 \to \ell^2} \geq a_0 > 4$,
we assume without loss of generality that $s\geq \frac{1}{4}$. Let $x = (x_n)$
be a sequence of real numbers defined by
\[
x_n = \left\{
\begin{array}{ll}
\qquad  1  & ~\mbox{ if }~ n = 0, \\  s (n+s) K_n & ~\mbox{ if }~ n \geq 1,
\end{array}
\right.
\]
where
\[
K_n = \frac{ \Gamma(\beta) \Gamma(n+\beta-\alpha) }{ \Gamma(n+\beta+1)
	\Gamma(\beta-\alpha+1) }
\]
and
\[
\alpha = \frac{2}{ 4+\varepsilon-\sqrt{(4+\varepsilon)\varepsilon} }, \quad
\beta = \frac{s^2}{ \alpha((4+\varepsilon)s-1) },
\]
for $\varepsilon > 0$. Using Stirling's formula, we see that $K_n \asymp
n^{-\alpha-1}$ and thus $x_n \asymp n^{-\alpha}$, since $x_n \asymp n K_n$.
Therefore, since
\[
\alpha = \frac{2}{ 4+\varepsilon-\sqrt{(4+\varepsilon)\varepsilon} } >
\frac{1}{2} \iff \varepsilon > 0,
\]
we have $x \in \ell^2$. Write $y = A_s x$. Hence,
\begin{align*}
y_n =&~ a_n \sum_{j=0}^{n} x_j + \sum_{j=n+1}^\infty a_j x_j\\
=&~ a_n + s \left( a_n \sum_{j=1}^{n}  (j+s) K_j + \sum_{j=n+1}^\infty K_j
\right)\\
=&~ a_n + s\frac{ \Gamma(\beta) }{ \Gamma(\beta-\alpha+1) } \left( a_n
\sum_{j=1}^{n} (j+s) \frac{ \Gamma(j+\beta-\alpha) }{ \Gamma(j+\beta+1) }
\right.\\
+&~ \left. \sum_{j=n+1}^\infty \frac{ \Gamma(j+\beta-\alpha) }{
	\Gamma(j+\beta+1) } \right).
\end{align*}
We know from a combinatorial identity that
\[
\sum_{j=1}^{n} \frac{ \Gamma(j+b) }{ \Gamma(j+c) }=\frac{ \Gamma(n+b+1) }{
	(1+b-c) \Gamma(n+c) }-\frac{ \Gamma(b+1) }{ (1+b-c) \Gamma(c) }.
\]
Thus, we can use this fact to find that
\[
y_n = (4+\varepsilon)s \frac{ \Gamma(n+\beta-\alpha+1) \Gamma(\beta) }{ (n+s)
	\Gamma(n+\beta) \Gamma(\beta-\alpha+1) }, \qquad (n\geq0).
\]
Observe that for $n=0$, we have $y_0 = 4+\varepsilon$ so $y_0 \geq
(4+\varepsilon) x_0$. We will now show that this inequality is true for all
$n\geq1$, as long that $\frac{1}{4} \leq s < s^*$. We can write
\[
y_n = (4+\varepsilon) \frac{ (n+\beta-\alpha) (n+\beta) }{ (n+s)^2 } x_n, \qquad
(n \geq 1).
\]
Thus we need to show that $(n+\beta-\alpha) (n+\beta) \geq (n+s)^2$, i.e., that
\begin{align}\label{eq5}
h_{\varepsilon,s}(n) := (2(\beta-s)-\alpha) n + \beta (\beta-\alpha) - s^2 \geq
0 , \qquad (n \geq 1).
\end{align}
Observe that $\alpha ,\, \beta$ are right-continuous function relative to
$\varepsilon$ at $\varepsilon = 0$. Hence, if we define $g(\varepsilon) :=
2(\beta-s)-\alpha$, we have
\[
|g(\varepsilon) - g(0)| < \eta
\]
provided that $0<\varepsilon<\mu$, where $\mu=\mu(\eta)>0$. However, $g(0) =
\frac{ 1-8s^2 }{ 2(4s-1) } > 0$ if $\frac{1}{4} < s < s^*$; thus we can set
$\eta = g(0)$ to be assured that there exist an $\varepsilon$ small enough so
that $g(\varepsilon) > 0$.

We have just shown that the leading coefficient of \eqref{eq5} is non-negative. Therefore, we just have to
make sure that $h_{\varepsilon,s}(1) \geq 0$ to make sure that \eqref{eq5} hold. Similarly to what we just did, we
can write $f(\varepsilon) := h_{\varepsilon,s}(1)$ and observe that $f$ is a
right continuous function relative to $\varepsilon$ at $\varepsilon = 0$. So if
$f(0) > 0$, we are assured that there exist a small enough $\varepsilon$ such
that $f(\varepsilon) = h_{\varepsilon,s}(1) \geq 0$. A computation gives us
\[
f(0) = \frac{ -24s^4-24s^3+8s^2+4s-1 }{ 2(4s-1)^2 }.
\]
Now, a simple analysis of this equation shows that $f(0) > 0$ for $\frac{1}{4} \leq s \leq
s^*$, with equality if and only if $s = s^*$. Hence, $y_n \geq
(4+\varepsilon) x_n$ for every $n \geq 0$, provided that $s<s^*$ and $\varepsilon$ is small enough.

Since $y = A_s x$ is a positive sequence and $x \in \ell^2$,
\[
\|A_s x\|_2^2 = \sum_{n=0}^\infty y_n^2 \geq (4+\varepsilon)^2 \sum_{n=0}^\infty
x_n^2 = (4+\varepsilon)^2 \|x\|_2^2,
\]
if $\varepsilon$ is small enough and $s<s^*$. It follows that $\|A\|_{\ell^2 \to \ell^2} > 4$ for all $s\in (0 ,\, s^*)$ and
thus, $s_0 \geq s^*$.

\section{The $p$-norm of $A_s$}\label{S:new-result}
In this section, we present the proof of Theorem \ref{T:As-p}. Parallel to the definition of $A_s$, consider the generalized Ces\`aro matrix
\[
C_s =
\begin{pmatrix}
\frac{1}{s}\, & 0\, & 0\, & 0\, & \cdots \\[3pt]
\frac{1}{1+s}\, & \frac{1}{1+s}\, & 0\, & 0\, & \cdots \\[3pt]
\frac{1}{2+s}\, & \frac{1}{2+s}\, & \frac{1}{2+s}\, & 0\, & \cdots \\[3pt]
\frac{1}{3+s}\, & \frac{1}{3+s}\, & \frac{1}{3+s}\, & \frac{1}{3+s}\, & \cdots
\\
\vdots\, & \vdots\, & \vdots\, & \vdots\, & \ddots \\
\end{pmatrix},
\]
and the generalized Copson matrix $C_s^{tr}$. We know from \cite{MR1933023} that, for $s\geq1$, $\|C_s\|_{\ell^p \to \ell^p}=q.$ It also follows that
\[
\|C_s^{tr}\|_{\ell^p \to \ell^p} = \|C_s\|_{\ell^q \to \ell^q} = p
\]
where $q$ is the Hölder conjugate of $p$, i.e., $1/p+1/q=1$.

As the first step, note that each of the entries of $A_s$ are positive and less
than or equal to those of $ C_s + C_s^{tr} $. Hence, we have
\begin{eqnarray*}
	\|A_s\|_{\ell^p \to \ell^p} &\leq&  \|C_s + C_s^{tr}\|_{\ell^p \to \ell^p} \\
	&\leq& \|C_s\|_{\ell^p \to \ell^p} + \|C_s^{tr}\|_{\ell^p \to \ell^p} \\
	&=& q + p = p q = \frac{p^2}{p-1}.
\end{eqnarray*}
We then proceed to show that $\|A_s\|_{\ell^p \to \ell^p} \geq \tfrac{p^2}{p-1}.$ Let
\[
x_m := \left( s^{-\frac{1}{p}},\, (1+s)^{-\frac{1}{p}},\,
(2+s)^{-\frac{1}{p}},\, \cdots,(m+s)^{-\frac{1}{p}},\, 0,\,0,\, \cdots
\right)^{tr},
\]
We then have
\begin{align*}
\|A_s x_m\|_{\ell^p \to \ell^p}^p
&\geq \sum_{n=0}^m \left( \frac{1}{n+s} \sum_{k=0}^{n-1} (k+s)^{-\frac{1}{p}} +
\sum_{k=n}^m (k+s)^{-1-\frac{1}{p}} \right)^p\\
&\geq \sum_{n=0}^m \left( \frac{1}{n+s} \int_{0}^{n} (x+s)^{-\frac{1}{p}} dx +
\int_{n}^{m+1} (x+s)^{-1-\frac{1}{p}} dx \right)^p\\
&= (pq)^p \sum_{n=0}^m \frac{1}{n+s} \left(  1-\frac{1}{p} \Big(\frac{s}{n+s}
\Big)^\frac{1}{q} - \frac{1}{q} \Big(\frac{m+s+1}{n+s} \Big)^{-\frac{1}{p}}
\right)^p.
\end{align*}
Since the summand in the second line of the chain of equations is positive,
$\left(  1-\frac{1}{p} \left(\frac{s}{n+s} \right)^{1/q} - \frac{1}{q}
\left(\frac{m+s+1}{n+s} \right)^{-1/p} \right)$ must also be positive. Hence, we
can use Bernouilli inequality to deduce that $\|A_s x_m\|_{\ell^p \to \ell^p}^p$
is
\begin{align*}
&\geq (pq)^p \sum_{n=0}^m \frac{1}{n+s} \left( 1-p\left( \frac{1}{p}
\Big(\frac{s}{n+s} \Big)^\frac{1}{q} + \frac{1}{q} \Big(\frac{m+s+1}{n+s}
\Big)^{-\frac{1}{p}} \right) \right)\\
&= (pq)^p \|x_m\|_{\ell^p \to \ell^p}^p - \gamma_m,
\end{align*}
where
\[
\gamma_m := (pq)^p \sum_{n=0}^m \frac{1}{n+s} \left(  \Big(\frac{s}{n+s}
\Big)^\frac{1}{q} + \frac{p}{q} \Big(\frac{m+s+1}{n+s} \Big)^{-\frac{1}{p}}
\right).
\]
This implies that
\[
\|A_s\|_{\ell^p \to \ell^p}^p \geq \frac{\|A_s x_m\|_{\ell^p \to
		\ell^p}^p}{\|x_m\|_{\ell^p \to \ell^p}^p} \geq (pq)^p -\frac{
	\gamma_m}{\|x_m\|_{\ell^p \to \ell^p}^p}.
\]
It is enough now show that $\gamma_m/\|x_m\|_{\ell^p \to \ell^p}^p \to 0$
whenever $m\to \infty$. First, note that
\begin{align}\label{eq_2}
\lim_{m\to \infty} \|x_m\|_{\ell^p \to \ell^p}^p = \lim_{m\to \infty}
\sum_{n=0}^m \frac{1}{n+s} = \infty.
\end{align}
Moreover, we have
\begin{align*}
\gamma_m &= (pq)^p \sum_{n=0}^m \left(  \frac{1}{n+s} \Big(\frac{s}{n+s}
\Big)^\frac{1}{q} + \frac{p}{q} \frac{1}{n+s} \Big(\frac{m+s+1}{n+s}
\Big)^{-\frac{1}{p}} \right)\\
&\leq c_1 \sum_{n=0}^m (n+s)^{-1/q-1} + c_2 (m+s+1)^{-1/p} \sum_{n=0}^m
(n+s)^{1/p-1}\\
&\leq c_1' (m+s+1)^{-1/q} + c_2' (m+s+1)^{-1/p} \cdot (m+s+1)^{1/p}\\[1pt]
&\leq c_1' (s+1)^{-1/q} + c_2'.
\end{align*}
Thus, the sequence $\gamma_m$ is bounded and there exist a constant $c$ such
that $\gamma_m \leq c$ for every $m\geq0$. Hence,
\[
0 \leq \frac{\gamma_m}{\|x_m\|_{\ell^p \to \ell^p}^p} \leq
\frac{c}{\|x_m\|_{\ell^p \to \ell^p}^p}.
\]
From \eqref{eq_2}, it follows that $\gamma_m/\|x_m\|_{\ell^p \to \ell^p}^p \to 0$
whenever $m\to \infty$. Therefore, $\|A_s\|_{\ell^p \to \ell^p} \geq pq = \tfrac{p^2}{p-1}$
and we are done.

\section{The norm of a special lacunary C-matrix} \label{S:cond-an12}
In this section, we present the proof of Theorem \ref{C:lacunary-Lmatrix-2}. Suppose that $y=Cx$. Then, we have
\[
N^{n/2} y_{N^n} = \sum_{j=0}^{N} x_j + \sum_{j=N+1}^{N^2} x_j + \cdots +
\sum_{j=N^{n-1}+1}^{N^{n}} x_j
\]
and $y_k=0$ for the other indices $k$. By Cauchy--Schwartz,
\begin{align*}
N^{n/2} y_{N^n} \leq&~ (N+1)^{\frac{1}{2}} \left( \sum_{j=0}^{N} x_j^2
\right)^{\frac{1}{2}} + (N^{2}-N)^{\frac{1}{2}} \left(\sum_{j=N+1}^{N^2} x_j^{2}
\right)^{\frac{1}{2}} \\
&+ \cdots + (N^{n}-N^{n-1})^{\frac{1}{2}} \left(\sum_{j=N^{n-1}+1}^{N^{n}}
x_j^{2} \right)^{\frac{1}{2}}.
\end{align*}
Once more, use the Cauchy--Schwartz inequality to get
\begin{align*}
N^{n/2} y_{N^n} \leq&~ \Big[ (N+1)^t + (N^2-N)^{\frac{1}{2}} + \cdots +
(N^{n}-N^{n-1})^{\frac{1}{2}} \Big]^{\frac{1}{2}} \cdot \\
& ~\Big[(N+1)^{1-t} \sum_{j=0}^{N} x_j^2  + (N^{2}-N)^{\frac{1}{2}}
\sum_{j=N+1}^{N^2} x_j^{2} \\
& \quad+ \cdots + (N^{n}-N^{n-1})^{\frac{1}{2}} \sum_{j=N^{n-1}+1}^{N^{n}}
x_j^{2} \Big]^{\frac{1}{2}}.
\end{align*}
Write
\begin{align*}
B_n:=~& (N+1)^t + (N^2-N)^{\frac{1}{2}} + \cdots +
(N^{n}-N^{n-1})^{\frac{1}{2}}\\
=~& (N+1)^t + \sqrt{N-1} \, \frac{\sqrt{N}^{\,n}-\sqrt{N}}{\sqrt{N}-1}.
\end{align*}
Then, for each $0 \leq t \leq 1$, we have
\begin{align*}
y_{N^n}^{2} ~\leq~& \frac{B_n(N+1)^{1-t}}{N^{n}} \sum_{j=0}^{N} x_j^2 +
\frac{B_n \sqrt{N^2-N}}{N^{n}} \sum_{j=N+1}^{N^2} x_j^2  \\
&~+ \cdots + \frac{B_n \sqrt{N^{n}-N^{n-1}}}{N^{n}} \sum_{j=N^{n-1}+1}^{N^{n}}
x_j^{2}.
\end{align*}
Therefore,
\begin{align}\label{eq: up_bound_B}
\|Cx\|_2^2 &=  \sum_{n=0}^{\infty} y_n^2
= \sum_{n=1}^{\infty} y_{N^n}^2\notag\\
&\leq \eta_0 \sum_{j=0}^{N} x_j^2 + \eta_1 \sum_{j=N+1}^{N^2} x_j^2 + \eta_2
\sum_{j=N^2+1}^{N^3} x_j^2 + \cdots ,
\end{align}
where
\begin{align*}
\eta_0 = \sum_{n=1}^\infty \frac{B_n(N+1)^{1-t}}{N^n}
&= \sum_{n=1}^\infty \frac{1}{N^n} \left( (N+1)^t + \sqrt{N-1} \,
\frac{\sqrt{N}-\sqrt{N}^{\,n}}{1-\sqrt{N}} \right)\\
&= (N+1)^{1-t}\frac{(N+1)^t(\sqrt{N}-1)+\sqrt{N-1}}{(\sqrt{N}-1)(N-1)},
\end{align*}
and
\begin{align*}
\eta_k &= \sum_{n=k+1}^\infty \frac{B_n\sqrt{N^{k+1}-N^k}}{N^n}
=\sqrt{N-1}\sqrt{N}^{\,k} \sum_{n=k+1}^\infty \frac{B_n}{N^n} \\
&=\frac{N-1}{(\sqrt{N}-1)^{2}}+\left(\frac{(N+1)^{t}}{\sqrt{N-1}}-\frac{\sqrt{N}}{\sqrt{N}-1}\right)\frac{1}{\sqrt{N}^{\,k}},
\qquad (k \geq 1).
\end{align*}
The upper bound \eqref{eq: up_bound_B} is valid for any $t\in[0,1]$. However,
the optimal $t$ is
\[
t=1-\log_{N+1}\sqrt{N-1}
\]
for which $\eta_k$ reduces to
\begin{align*}
\eta_k = \frac{N-1}{(\sqrt{N}-1)^{2}}-\frac{1}{(\sqrt{N}+1)\sqrt{N}^{\,k}},
\qquad (k \geq 0).
\end{align*}
Thus, we have
\begin{align*}
\eta_k  \leq \frac{N-1}{(\sqrt{N}-1)^{2}}, \qquad (k \geq 0).
\end{align*}
This special choice of $t$ implies
\begin{align*}
\|Cx\|_2^2 \leq \frac{N-1}{(\sqrt{N}-1)^2} \left( \sum_{j=0}^{N} x_j^2 +
\sum_{j=N+1}^{N^2} x_j^2 + \cdots \right) = \frac{N-1}{(\sqrt{N}-1)^2} \|x\|_2^2.
\end{align*}
Hence,
\[
\|C\|_{\ell^2 \to \ell^2} \leq \frac{\sqrt{N-1}}{\sqrt{N}-1}.
\]
We now proceed to show that this upper bound is attained.

Let $x:=(x_k)$ be the following sequence: $x_k$ equals $1$ for indices between
$0$ and $N$, and equals $1/\sqrt{N}^{\,n}$ for every indices between $N^{n}+1$
to $N^{n+1}$, for $n=1,2,\dots,m$, and finally equals $0$ everywhere else. Then
\begin{align*}
\|x\|_2^2 &=  \sum_{j=0}^{N} x_j^2 + \sum_{j=N+1}^{N^2} x_j^2 + \cdots +
\sum_{j=N^{m-1}+1}^{N^{m}} x_j^2\\
&= N+1 + (N^2-N)\frac{1}{N} + \cdots + (N^m-N^{m-1})\frac{1}{N^{m-1}}\\
&= 2+(N-1)m.
\end{align*}
Once more, write $y=Cx$. Then
\begin{align*}
N^{n/2} y_{N^n} &= \sum_{j=0}^{N} x_j + \sum_{j=N+1}^{N^2} x_j + \cdots +
\sum_{j=N^{n-1}+1}^{N^{n}} x_j\\
&=N+1 + (N^2-N)\frac{1}{\sqrt{N}} + (N^3-N^2)\frac{1}{\sqrt{N}^{\,2}} \\
&\quad+ \cdots + (N^n-N^{n-1})\frac{1}{\sqrt{N}^{\,n-1}}\\
&= (\sqrt{N}+1)\sqrt{N}^{\,n} - (\sqrt{N}-1),
\end{align*}
for every $n\leq m$, and thus
\begin{align*}
y_{N^n}^2 = \left( (\sqrt{N}+1) - \frac{\sqrt{N}-1}{\sqrt{N}^{\,n}} \right)^{2}.
\end{align*}
Hence,
\begin{align*}
\|Cx\|_2^2 &\geq \sum_{n=1}^m y_{N^n}^2 = \sum_{n=1}^m \left( (\sqrt{N}+1) -
\frac{\sqrt{N}-1}{\sqrt{N}^{\,n}} \right)^{2}\\
&\geq (\sqrt{N}+1)^2 m + c,
\end{align*}
for a certain constant $c$ (for example, $c=-2(\sqrt{N}+1)$ works). Therefore,
\begin{align*}
\|C\|_{\ell^2 \to \ell^2}^2 \geq \frac{\|Cx\|_2^2}{\|x\|_2^2} \geq
\frac{(\sqrt{N}+1)^2 m + c}{2+(N-1)m} 
\end{align*}
for every $m\geq 1$. By letting $m \to \infty$, we get
\begin{align*}
\|C\|_{\ell^2 \to \ell^2}^2 \geq \frac{(\sqrt{N}+1)^2}{N-1} =
\frac{N-1}{(\sqrt{N}-1)^2}.
\end{align*}
Therefore,
\[
\|C\|_{\ell^2 \to \ell^2} \geq \frac{\sqrt{N-1}}{\sqrt{N}-1},
\]
and thus the equality follows.

\section{Concluding remarks}
\begin{enumerate}[(i)]
\item The precise value of $s_0$ is still unknown. Find $s_0$?
\item In the light of Theorem \ref{T:As-p}, we define
\[
s_0^{(p)} := \inf \{s: \|A_s\|_{\ell^p \to \ell^p} = \tfrac{p^2}{p-1} \}.
\]
That theorem ensures
\[
s_0^{(p)} \leq 1.
\]
Find $s_0^{(p)}$.
\item Find $\|A_s\|_{\ell^p \to \ell^q}$, where $p \ne q$.
\item As in the case $p=q$, does $\|A_s\|_{\ell^p \to \ell^q}$ remain constant for large values of $s$? If so, we define $s_0^{(pq)}$ by slightly modifying the definition of $s_0^p$ given in (ii). Then how does  $s_0^{(pq)}$ depend on the parameters $p$ and $q$?
\end{enumerate}

\section*{Declaration of competing interest}
The authors declare that they have no competing interests.

\section*{Acknowlegdment}
For this work, the first author received the USRA research award. The
second author was supported by the NSERC Discovery Grant (Canada).
	

\bibliographystyle{elsarticle-num}
\bibliography{The-critical-point-references}

\begin{thebibliography}{10}
\expandafter\ifx\csname url\endcsname\relax
  \def\url#1{\texttt{#1}}\fi
\expandafter\ifx\csname urlprefix\endcsname\relax\def\urlprefix{URL }\fi
\expandafter\ifx\csname href\endcsname\relax
  \def\href#1#2{#2} \def\path#1{#1}\fi

\bibitem{MRtmp}
L.~Bouthat, J.~Mashreghi, The norm of an infinite {L}-matrix, Oper. Matrices
  15~(1) (2021) 47--58.
\newblock \href {http://dx.doi.org/10.7153/oam-2021-15-04}
  {\path{doi:10.7153/oam-2021-15-04}}.

\bibitem{MR4032210}
J.~Mashreghi, T.~Ransford, Hadamard multipliers on weighted {D}irichlet spaces,
  Integral Equations Operator Theory 91~(6) (2019) Paper No. 52, 13.
\newblock \href {http://dx.doi.org/10.1007/s00020-019-2551-1}
  {\path{doi:10.1007/s00020-019-2551-1}}.

\bibitem{MR2500010}
J.~Mashreghi, Representation theorems in {H}ardy spaces, Vol.~74 of London
  Mathematical Society Student Texts, Cambridge University Press, Cambridge,
  2009.
\newblock \href {http://dx.doi.org/10.1017/CBO9780511814525}
  {\path{doi:10.1017/CBO9780511814525}}.

\bibitem{MRtmp2}
L.~Bouthat, J.~Mashreghi, {L}-matrices with lacunary coefficients, Oper.
  Matrices. to appear.

\bibitem{MR1554854}
D.~Hilbert, Ein {B}eitrag zur {T}heorie des {L}egendre'schen {P}olynoms, Acta
  Math. 18~(1) (1894) 155--159.
\newblock \href {http://dx.doi.org/10.1007/BF02418278}
  {\path{doi:10.1007/BF02418278}}.

\bibitem{MR0030620}
G.~H. Hardy, Divergent {S}eries, Oxford, at the Clarendon Press, 1949.

\bibitem{MR0305568}
D.~M. Young, Iterative solution of large linear systems, Academic Press, New
  York-London, 1971.

\bibitem{MR701570}
M.~D. Choi, Tricks or treats with the {H}ilbert matrix, Amer. Math. Monthly
  90~(5) (1983) 301--312.
\newblock \href {http://dx.doi.org/10.2307/2975779}
  {\path{doi:10.2307/2975779}}.

\bibitem{MR998639}
H.~C. Rhaly, Jr., Terraced matrices, Bull. London Math. Soc. 21~(4) (1989)
  399--406.
\newblock \href {http://dx.doi.org/10.1112/blms/21.4.399}
  {\path{doi:10.1112/blms/21.4.399}}.

\bibitem{MR3112183}
H.~C. Rhaly, Jr., Hyponormality-preserving finite rank perturbations of
  terraced matrices, J. Nigerian Math. Soc. 32 (2013) 281--288.

\bibitem{MR1034416}
B.~E. Rhoades, Lower bounds for some matrices. {II}, Linear and Multilinear
  Algebra 26~(1-2) (1990) 49--58.
\newblock \href {http://dx.doi.org/10.1080/03081089008817965}
  {\path{doi:10.1080/03081089008817965}}.

\bibitem{MR3458973}
I.~Almasri, Absolutely summing terraced matrices, Concr. Oper. 3~(1) (2016)
  1--7.
\newblock \href {http://dx.doi.org/10.1515/conop-2016-0001}
  {\path{doi:10.1515/conop-2016-0001}}.

\bibitem{MR3527878}
N.~Durna, M.~Yildirim, Generalized terraced matrices, Miskolc Math. Notes
  17~(1) (2016) 201--208.
\newblock \href {http://dx.doi.org/10.18514/MMN.2016.1272}
  {\path{doi:10.18514/MMN.2016.1272}}.

\bibitem{MR4033299}
M.~\.{I}lkhan, Norms and lower bounds of some matrix operators on {F}ibonacci
  weighted difference sequence space, Math. Methods Appl. Sci. 42~(16) (2019)
  5143--5153.
\newblock \href {http://dx.doi.org/10.1002/mma.5244}
  {\path{doi:10.1002/mma.5244}}.

\bibitem{MR3641772}
M.~Jevti\'{c}, B.~Karapetrovi\'{c}, Hilbert matrix on spaces of {B}ergman-type,
  J. Math. Anal. Appl. 453~(1) (2017) 241--254.
\newblock \href {http://dx.doi.org/10.1016/j.jmaa.2017.04.002}
  {\path{doi:10.1016/j.jmaa.2017.04.002}}.

\bibitem{MR3545935}
O.~F. Brevig, K.-M. Perfekt, K.~Seip, A.~G. Siskakis, D.~Vukoti\'{c}, The
  multiplicative {H}ilbert matrix, Adv. Math. 302 (2016) 410--432.
\newblock \href {http://dx.doi.org/10.1016/j.aim.2016.07.019}
  {\path{doi:10.1016/j.aim.2016.07.019}}.

\bibitem{MR3208804}
I.~Chalendar, J.~R. Partington, Norm estimates for weighted composition
  operators on spaces of holomorphic functions, Complex Anal. Oper. Theory
  8~(5) (2014) 1087--1095.
\newblock \href {http://dx.doi.org/10.1007/s11785-013-0340-4}
  {\path{doi:10.1007/s11785-013-0340-4}}.

\bibitem{MR1933023}
C.-P. Chen, D.-C. Luor, Z.-y. Ou, Extensions of {H}ardy inequality, J. Math.
  Anal. Appl. 273~(1) (2002) 160--171.
\newblock \href {http://dx.doi.org/10.1016/S0022-247X(02)00232-9}
  {\path{doi:10.1016/S0022-247X(02)00232-9}}.

\end{thebibliography}
\end{document}